\documentclass[11pt]{article}
\usepackage[T1]{fontenc}
\usepackage[utf8]{inputenc}

\usepackage{amsmath, amssymb}

\usepackage{amsthm}

\usepackage[a4paper, margin=1in]{geometry}
\usepackage{setspace}
\usepackage{minted}
\usepackage{newpxtext}
\usepackage{authblk}
\usepackage[hyperindex,breaklinks]{hyperref}
\usepackage{xcolor}
\usepackage{cleveref}
\usepackage{hyperref}
\usepackage{footmisc}

\usepackage{fancyvrb}
\usepackage{fontspec}
\usepackage{listings}
\usepackage{color}

\definecolor{keywordcolor}{rgb}{0.7, 0.1, 0.1}   
\definecolor{tacticcolor}{rgb}{0.0, 0.1, 0.6}    
\definecolor{commentcolor}{rgb}{0.4, 0.4, 0.4}   
\definecolor{symbolcolor}{rgb}{0.0, 0.1, 0.6}    
\definecolor{sortcolor}{rgb}{0.1, 0.5, 0.1}      
\definecolor{attributecolor}{rgb}{0.7, 0.1, 0.1} 

\title{A correspondence problem for mathematical proof}

\author[1,2]{Simon DeDeo\textsuperscript{$\Delta$}\thanks{Email: sdedeo@andrew.cmu.edu}}
\author[3,4]{Eamon Duede\textsuperscript{$\Delta$}\thanks{Email: eduede@purdue.edu}}
\affil[1]{Carnegie Mellon University}
\affil[2]{Santa Fe Institute}
\affil[3]{Purdue University}
\affil[4]{Argonne National Laboratory}

\date{} 

\makeatletter
\newcommand{\setword}[2]{%
  \phantomsection
  #1\def\@currentlabel{\unexpanded{#1}}\label{#2}%
}
\makeatother

\begin{document}
\linespread{1.25}
\maketitle
{
\renewcommand{\thefootnote}{$\Delta$}
\footnotetext{These authors contributed equally to this work.}
}

\maketitle

\begin{abstract}
\noindent Mathematical proofs are often said to justify their conclusions by indicating the existence of a corresponding formal derivation. We argue that this widespread view relies on an under-examined notion of correspondence, or what it means for a particular derivation to ``correspond'' to a particular proof. Mere existence of a formalization is not enough, and a substantive account of the required correspondence resolves into two criteria---adequate representation (of the original theorem) and tracking (of the steps in the original proof). An examination of the actually-existing formalization systems we have today shows the variety of quasi-empirical ways we establish these criteria, and points towards new burdens that may be placed on the future evolution of mathematics itself.
\end{abstract}

\section*{Preliminaries}
\noindent In general, it is uncontroversial to claim that mathematical proofs possess epistemic value. Controversy arises, however, as soon as one attempts to articulate what that value consists in, and how proofs manage to secure it. After all, proofs plainly do many things: communicate mathematical ideas, transmit understanding, and unify or extend previously distinct bodies of knowledge. Yet, proofs are also taken to play distinctive justificatory roles. In this mode, proofs are supposed to establish the \textit{correctness} of mathematical claims and thereby warrant the striking degree of consensus characteristic of mathematical practice. What has proved philosophically puzzling is that the mathematical proofs one finds in top journals and textbooks routinely succeed in this role despite being informal, incomplete, saturated with tacit assumptions, gaps, appeals to background understanding, and not infrequently, peppered with local mistakes. How such objects can nonetheless justify mathematical claims with the apparent reliability and deductive security of mathematical knowledge is far from obvious.

The dominant response has been to treat the justificatory force of an informal mathematical proof as grounded, ultimately, in the existence of a corresponding formal derivation. According to this so-called ``Standard View'', an informal proof is correct insofar as it indicates that a fully explicit proof, in a suitable formal proof system, could be produced \cite{azzouni2004derivation, avigad2019automated, hamami2022mathematical, mac1986functions}. The Standard View comes in many flavors, but what they all share is a commitment to the idea that beneath the informal presentation of a \textit{genuine proof}, there exists, or ought to exist, an ``unbroken chain of logical inferences from an explicit set of axioms.'' \cite[p.x]{hales2012dense}. Opponents of the Standard View generally deny that the correctness of informal proofs can be grounded simply in the existence of an underlying derivation, emphasizing instead the irreducibly semantic, methodological, and practice-bound features of proof that are not captured by a corresponding unbroken chain of formal inferences \cite{rav1999we, tanswell2015problem}.

What both proponents and critics of the Standard View share, though, is a largely tacit commitment to the antecedent idea that, for any given proof, it makes sense to speak of a formal proof that \textit{corresponds} to an informal counterpart. Given this shared commitment, the debate has then turned on whether, to what extent, or how, a that corresponding formal derivation contributes to the epistemic value of justification.

In this paper, we argue that the antecedent idea of correspondence itself cannot be taken for granted. The issue, which is of more than theoretical interest for contemporary library-based formalization projects, is not whether, for a given, informally expressible theorem, a formal proof exists, but, rather, what it could mean for such a proof to correspond to a \textit{particular} human proof-idea at all. That is, successfully deriving a formal proof which corresponds not merely to a proof of some \textit{particular theorem}, but, rather to some \textit{particular proof} of that theorem is itself a kind of epistemic achievement, one that contemporary formalization ecosystems may systematically frustrate.

Our paper is divided into five sections. First, in Section~\ref{correspondence}, we briefly rehearse and situate the Standard View and its critics, emphasizing their shared commitment to the intelligibility of a correspondence relation between informal and formal proofs. We distinguish between a ``thin'' and ``thick'' notion of correspondence, and argue for the primacy of the thick in how correspondence has appeared in the literature to date. In Section~\ref{metaproblem}, we unfold this thick notion: for a correspondence to do what is required, the formal proof must both (1) represent the claims in an adequate fashion, and (2) track the steps of the informal proof of those claims. Both are challenging tasks, but we discuss how actually-existing formalization projects shift the burden of correspondence to the representational side as much as possible, and why this is a good idea. Next, in Section~\ref{leanturn}, we show how inherited infrastructure and entrenched dependency structures further destabilize the traditional picture of correspondence by enabling formally correct proofs that hybridize fragments from distinct modeling traditions. We argue that this yields a deeper asymmetry. While human proofs are semantically upstream of formal derivations, formal derivations assembled within large libraries may, and likely do, detach from the the proof-ideas and explanatory structures they might be taken to express. 

In the light of this we present a positive account, in Section~\ref{routes}, of three routes to believing a correspondence has been obtained: causal, explanatory, and structural. This places formalization on a more equal footing with other human proofs, as a practice that systematically constrains proof-space, selecting particular definitions and derivational pathways as intelligible, explanatory, and reusable relative to shared background models within a mathematical community. These constraints emerge in part from the experience of error-making and dialogue between mathematicians, and have a genealogy that responds just as much to the cognitive constraints and histories of the participants as to the underlying structure of proof space. We draw, briefly, on a recent controversy in the formalization community, to show how, and for what reasons, practitioners navigate the correspondence problem. We conclude in Section~\ref{ai}.

\section{Correspondence}
\label{correspondence}

The proofs that appear in leading mathematical journals and textbooks are a mixture of natural language and mathematical symbolism, structured around a shared understanding that, out of necessity, leaves many inference steps implicit. Following the convention of the literature, we refer to proofs of this kind as \textit{informal proofs}. \textit{Formal} proofs, by contrast, are derivations carried out within a fixed proof system, in which every deductive dependency from axioms to conclusion is made precise and, in principle, mechanically checkable. That informal and formal proofs both aim to establish mathematical theorems suggests that, at least in principle, a given result may admit both an informal and a formal proof. Moreover, it has become a kind of mathematical orthodoxy that if an informal proof is genuinely correct, then a \textit{corresponding} formal derivation of the theorem must exist. This orthodoxy is visible in, for example, the legacy of the mid-century structuralist traditions of Mac~Lane and Bourbaki \cite{lane1935logical,bourbaki1970theorie}, and is likely continuous with the still earlier logicist and formalist ambitions of Frege \cite{frege1879begriffsschrift} and Hilbert \cite{hilbert1938grundzuge}, respectively.

This Standard View is not one thesis, but a family united by the appeal to formal derivation as the justificatory ground and desideratum for recognizing the epistemic value of a given informal proof. It is also vague, at least as rendered here, because it seems there are two ways to render the relationship between informal and formal proof. The first is to say that, if an informal proof $\pi$ of theorem $\theta$ is a genuine proof (and not merely an argument that \textit{seems}, for all the world, like a proof but isn't \cite{de2021groundwork}) of $\theta$, then there exists a formal derivation $\delta$ of $\theta$. Call this the \textit{thin reading}. The second is to say that if an informal proof $\pi$ of theorem $\theta$ is a genuine proof, then there exists a formal derivation $\delta$ of $\theta$ that corresponds to $\pi$. Call this the \textit{thick reading}. On the thin reading any formal derivation of $\theta$ will do. On the thick reading, correspondence demands that the derivation $\delta$ stand in some privileged structural, causal, or explanatory relation to the proof-idea represented in $\pi$ \textit{itself}. 

Proponents of the Standard View do not always mark this distinction, and the literature accordingly slides between the two readings in ways that obscure what is really at stake.\footnote{For instance, when Avigad says that ``with enough work [an informal proof] can be turned into a formal derivation...'' \cite{avigad2021reliability}, or Mac Lane says that a sufficiently detailed informal proof (which he calls a sketch) should ``make possible a routine translation of \textit{this} sketch into a formal proof'' \cite{mac1986functions}, both suggestions are clearly thicker than the mere existence of a formal derivation.} Surely, if an informal proof is correct, some formal derivation of $\theta$ exists. That's just to say that if $\theta$ is true (and provable), then $\theta$ is derivable in a suitable formal system equal in power to that of the informal case. The thin reading is trivially available but makes it hard to see what distinctive epistemic work the correctness of $\pi$ as some \textit{particular} proof of $\theta$ could be doing. Any other informal proof $\pi^\prime$ of $\theta$, however different in strategy, conceptual apparatus, or explanatory ambition, would be equally well served by the \textit{same} underlying derivation $\delta$. The thin reading, in other words, seems to validate theorems, not proofs---but it would seem that the correctness of specific proofs, rather than the truth, or provability, of specific theorems, is most commonly at issue~\cite{gowers2023makes}.

To see this, suppose we have two informal proofs, $\pi_1$ and $\pi_2$, of the same theorem $\theta$. The proof-idea and strategy in $\pi_1$ relies on a complex and cumbersome set of manipulations within a single area of mathematics. The proof-idea and strategy in $\pi_2$ relies on a simple and deft integration of seemingly disparate ideas from multiple areas of mathematics. Now, suppose mathematicians are fairly certain that $\pi_1$ is a genuine proof, but are much more excited about $\pi_2$ because of its elegance and unifying power. They want to know if $\pi_2$ is a genuine proof. So, they formalize it as $\delta$. On the thin reading, the formalization $\delta$ confirms only that $\theta$ is derivable. But, this is something mathematicians could already have inferred from the existence of $\pi_1$. Whether $\pi_2$'s distinctive strategy actually works, whether its elegant integration of disparate ideas constitutes a genuine proof, is left entirely unevaluated.

It would seem that, in order to give an account of the epistemic value of a particular proof as \textit{proof}, the Standard View requires the thick reading, and thus needs to say what correspondence relation needs to exist, between an informal proof $\pi$ and formal proof $\delta$, such that the formal structure of $\delta$ answers to $\pi$'s inferential organization. It needs to do so, moreover, in a way that would not equally answer to genuinely distinct proof $\pi^\prime$ of the same theorem---or even a flawed proof, $\pi^\star$. The philosophical burden shifts from showing that a derivation exists to explaining what it could mean for a derivation to correspond to a particular proof-idea at all.

\section{Two Parts of the Correspondence Problem}
\label{metaproblem}

For the most part, the philosophical literature has treated the notion of correspondence as an unanalyzed primitive, one sufficiently clear that it can be gestured at and appealed to without much scrutiny. Even prominent critics of the Standard View leave the antecedent intelligibility of correspondence intact. For instance, Rav's influential attack on the Standard View is motivated to undermine the idea that what is epistemically valuable about a mathematical proof $\pi$ can be \textit{exhausted} by its corresponding formal derivation $\delta$ \cite{rav1999we}. Even Hamami's careful reconstruction of the Standard View \cite{hamami2022mathematical} does not analyze what makes a formal derivation a derivation of a particular proof rather than merely a derivation of a particular theorem. The shared architecture of the debate, in other words, is one in which correspondence is a relation whose adequacy or extent may be disputed, but whose intelligibility is taken for granted.

The intuition behind the thick position is that for any particular informal proof $\pi$, the bare existence of a formal proof $\delta$ is not enough to convince us that $\pi$ deserves to be called a success. As we will argue, we hope in particular to convince ourselves that the statement established at the end of $\delta$ is an \emph{adequate representation} of the statement that appears at the beginning of $\pi$. We hope, too, to see how features of $\delta$ \emph{track} the reasoning of $\pi$. Both the informal and formal proof present their reasoning as an unfolding series of steps. For $\delta$ to track $\pi$ is for there to exist an appropriate map from each step in $\pi$ to some step, or set of steps, in $\delta$, one that preserves the dependency structure of the informal argument even where the formal proof might reorder or decompose individual steps. The idea of correspondence, then, is the joint satisfaction of tracking and representation conditions.

\begin{quote}
    \textbf{Correspondence}: a formal proof $\delta$ corresponds to an informal proof $\pi$ just in case it successfully tracks the inferential sequence argued for in $\pi$, and is an adequate representation of the theorem $\theta$ expressed by $\pi$.
\end{quote}

For tracking, the relevant issue is not securing line-by-line mimicry, but rather the preservation of enough inferential architecture that $\pi$ and $\delta$ can be recognized as implementing the same proof ideas and strategy. In recognizing the adequacy of representation, the relevant issue is not establishing that $\delta$ is a mere grammatical paraphrase of the theorem statement in $\pi$, but, rather the adequacy of the formal objects, predicates, and dependencies used to model the mathematical content at issue in $\pi$. 

Both parts of the correspondence problem have been simmering for decades. Early formalization projects, such as Bourbaki, were conducted by mathematicians at the very center of their fields, where it was an almost unstated assumption that the participants had adequate representations of the informal mathematics of their time, were able to contemplate their formalizations in a similar light, and were able to verify that the steps in the parallel proofs tracked in the appropriate fashion. The informal proof might expand as it was formalized (with, for instance, irritating asides to properly encode functions as sets) but if the formalization appeared to take significant detours, this would be a sign that the informal proof needed to be changed, or rewritten entirely.

Many of these otherwise unspoken needs for both representation and tracking were brought to wider attention by the advent of new computational tools for formal proof. Today, languages for formalization such as \texttt{Lean}, \texttt{Rocq}, and \texttt{Isabelle/HOL} are in active use by large communities and prominent mathematicians, and have formalized significant achievements in mathematics.\footnote{See, \emph{e.g.}, the Liquid Tensor Experiment~\cite{lte,hartnett2021proof} led by Johan Commelin, that formalized a sophisticated ``challenge'' theorem in Analytic Geometry by Fields Medalist Peter Scholze.}

To understand the two-part correspondence problem in the light of recent developments, we begin with a toy example. Consider this ``proof'' of Euclid's prime number theorem,
\begin{quote}
  \setword{\textbf{Theorem}}{primetheorem}: there is no largest prime number.

  \setword{\textbf{Proof$^\star$}}{proofdelta}:
  assume that there is a largest prime number, $N$. Construct the number $M$,
  equal to the product of all the numbers up to and including $N$, plus one.
  $M$ is not divisible by any number less than or equal to $N$. Therefore $M$
  is prime. $M$ is larger than $N$. Contradiction.
  \hfill $\square$
\end{quote}

Let us first consider challenges facing the \textit{tracking condition}. All of the major languages in use for formalization today are able to both state and prove Euclid's prime number \ref{primetheorem}. The existence of a (correct) formal proof of \ref{primetheorem} in, say \texttt{Lean}, however, tells us nothing about the informal proof articulated by \ref{proofdelta}. This is because \ref{proofdelta} is, of course, incorrect at sentence four: simply because $M$ is not divisible by any number less than or equal to $N$, that does not mean that $M$ is prime; the correct proof uses $M$ to establish the existence of a prime $p$, which need not be equal to $M$. As a result, the formal proof would not (assuming the system was sound) track the informal infelicity in sentence four and we would not find something that, in the formal language, amounted to a claim that $M$ is (always) prime.

Even when we have a \textit{correct} informal proof, however, the tracking problem remains. For any theorem, there are an indefinite number of formal proofs, which needn't all track the same informal proof. Euler's proof of the same theorem, for instance, uses the fact that the harmonic series, $\sum_{i=1}^\infty 1/n$, diverges, and with this clue as a guide, a human can use current formalization technology pick out a formal proof completely distinct from  Euclid's. There are also formal proofs that, at least on the surface, have no obvious, human-readable informal counterpart. Demonstrations of \texttt{AlphaProof} \cite{hubert2025olympiad} on Euclid's toy example produce valid \texttt{Lean} code, and thus valid proofs,\footnote{We take the validity of the \texttt{Lean} kernel, and the soundness of \texttt{mathlib}, as given, though this question is of course of broad interest \cite{mackenzie2004mechanizing}.} that, at least on the surface, has no clear tracking of the steps a human might take.

Things are no better with the \textit{adequate representation} condition. The natural language phrase ``there is no largest prime number'' in \ref{primetheorem} itself must be formalized. But we have no obvious criterion for the adequacy of the correspondence between this phrase and the supposed restatement of it in the formal language. In fact, in actually existing formalization ecosystems, when we type ``$\forall$ n, $\exists$ p, (Prime p) $\land$ (p > n)'', it is not always clear that the symbols and predicates really mean ``for all'', ``there exists'', ``is prime'' and ``greater than''.

These doubts are not merely a philosopher's fancy. An early attempt to formalize the law of quadratic reciprocity in \texttt{Isabelle} succeeded in tracking the informal proof, but failed due to representational inadequacy: the definition of ``prime'' was such that there were \textit{no primes}, so the proof amounted to the fact that everything is true of the empty set \cite{avigad2019automated}. This challenge persists today in what Johan Commelin has referred to as ``ceci n'est pas une pipe'' (see \cite{dedeo2024alephzero}). It is generally agreed that the solution is a quasi-empirical matter involving the ``stress-testing'' of definitions with a variety of test cases analogous to unit tests in computer engineering~\cite{garousi2016systematic}. In the words of one of the leaders of the Sphere Packing formalization project \cite{sphere},

\begin{quote}
    ``Stress-testing [definitions] is certainly not foolproof, and sometimes things don't get caught till you try something very bizarre or very deep. But if we balance careful design choices with clever testing then we can usually catch bugs.'' \cite{sid_quote}
\end{quote}

Evidence accumulates in favor of a correspondence up to a point of diminishing epistemic returns---rather than establishing it as a deductive certainty. The ``thick'' Standard View promises to underwrite a belief in the deductive certainty of an informal proof $\pi$, but we can gain only inductive belief of whether or not its promised justification of $\theta$ obtains.

In practice, actually-existing formalization projects seek to overcome the correspondence problem in two steps\footnote{We do not mean to assert this to be an acknowledged agenda for members of these communities.} which are naturally separated out in the ``type theoretic'' languages that dominate contemporary formalization. Consider a fragment of Euclid's proof of \ref{primetheorem} in \texttt{Lean},


\begin{lstlisting}
theorem primes_unbounded (n : Nat) : ∃ p, Prime p ∧ p ≥ n := by
  let m := Nat.factorial n + 1
  have h2 : 2 ≤ m := by
    have : 1 ≤ Nat.factorial n :=
      Nat.succ_le_of_lt (Nat.factorial_pos n) ...
\end{lstlisting}


The first line includes the ``type signature'': informally, {\tt primes\_unbounded} is defined as a function that takes a number, $n$ (of type {\tt Nat}, a natural number), and returns a proof (a term of type ``$\exists~p,~\textrm{Prime}~p\land p\geq n$'', which is, in turn, of type {\tt Prop}) of the existence of a number $p$ (of type {\tt Nat}, implicitly) with the property of being greater than or equal to $n$, and prime. Representational adequacy is then established by examination of the type signature. As long as one trusts the types, the signatures are usually easily verbalized; they can be thought of as ``declarations of intent'' in the sense that, at the end of the proof, the system itself verifies that a term of the promised type has, in fact, been constructed and returned.

The remaining lines, following {\tt by}, constitute the computational structure that actually constructs what the signature promises. Ideally, tracking adequacy is established by examining what follows {\tt by} to see if it mirrors one's own reasoning. In practice, however, when theorems become more complex, this becomes difficult. \texttt{Lean}, as with most type-theoretic languages, relies heavily ``tactics'': general-purpose tools that analyze what has come before their invocation and search within a restricted space to advance the proof to the next stage. The output of these tactics, as they unfold, can be inspected, but is often unintelligible and one relies, instead, on watching how, as a consequence of the tactic, the proof advances a little further towards the goal state. Proper \texttt{Lean} style suggests using the least powerful tactic necessary to move the proof to the desired next stage, but this is a stylistic choice, and nothing in the system prevents, for example, the use of an overpowered tactic which accomplishes the required steps in fashion that humans find genuinely mysterious.

The relative opacity of tracking compared to representation has real world consequences for how humans formalize. Consider, for example, the formalization of the polynomial Freiman–Ruzsa conjecture (PFR), whose informal proof \cite{gowers2025conjecture} runs approximately 13,000 words and contains seventeen statements labeled as theorem, corollary, lemma, or proposition. The corresponding \texttt{Lean} formalization \cite{pfr}, by contrast, has at least 910 corresponding `statements' known as \textit{declarations}. These are fragments of code that associate a type signature with a theorem, lemma, or other definition. Even if we broaden the criteria for `declaration' in the informal proof to include 67 numbered formulae, and perhaps an additional fifteen to twenty prose definitions, the formal proof is \textit{far more} declaration heavy, by a factor of at least ten.

The vast gap between the two modes is apparent. Actually-existing formalization projects have far more declarations of intent than their informal counterparts---which means that correspondence questions can lean more heavily on judgments of representational adequacy. This is, in turn, a somewhat happy state of affairs the nature of current type theoretic languages means that representational adequacy is often far easier to judge than tracking, and the correspondence becomes more assessable by exposing many intermediate claims in a readable form.

The large ratio is not unique to PFR. The same comparison using the formalization~\cite{sphere} of Maryna Viazovska's solution to the Sphere Packing Problem in Dimension 8 \cite{viazovska2017sphere}, finds 12 theorems and propositions, 65 numbered formulae, and perhaps twenty prose definitions in the informal proof compared to approximately 1362 declarations that appear in the human-authored formalization attempt, just before its unexpected completion by the A.I. tools of \emph{Math, Inc.} in late February 2026. 

The fully-unfolded formal proof will not necessarily label, or even preserve these boundaries. When compiled into its final form, for example, the resulting object is one large, directed graph, with potentially hundreds of thousands of nodes \cite{viteri2022epistemic}. The difficulty of working in a formal language, therefore, means that when humans attempt the formalization task, they must, out of necessity, work from a blueprint that tracks the original argument. 

The constraint has happy consequences for establishing the correspondences required for the epistemic weight of the Standard View to be felt. Nevertheless, those correspondences remain matters of evidence-gathering. Whether a formal derivation corresponds to a particular informal proof is, in the end, an empirical question---one whose resolution the Standard View presupposes but does not itself provide.

\section{Libraries, Lock-In, and Shared Models}
\label{leanturn}

We have seen that the Standard View requires a thick notion of correspondence and that this notion, in turn, requires that formal proofs adequately represent and track their informal counterparts. For any genuine informal proof $\pi$, furthermore, the establishment of a correspondence with a formal proof $\delta$ is in part an empirical achievement. The deductive confirmation of $\pi$ that the Standard View promises becomes, thereby, an inductive contingency, and formalizers have meaningful methods for securing it.

Even when formalizers are doing everything right, however, the formalization infrastructure they are working within can distort the correspondence by imposing proof-structural choices on authors that are artifacts of the library's history rather than features of $\pi$'s underlying proof-idea. Because (typical) informal proofs begin at a meaningful deductive distance from the axioms, they presuppose a large body of prior results that serve as an interface for what they wish to do next. The development of mathematics in its sociological and historical contingency is, therefore, constitutive of an ecosystem today that, in turn, evolves as new proofs are made.

In practice, when mathematicians turn to formalize an informal proof $\pi$, they do so not just within the intellectual ecosystem from which $\pi$ emerged, and in which $\pi$ is, hopefully, comprehensible and generative, but at the same time within the distinct ecosystem for $\delta$, which is a product of its own evolving historical and sociological contingencies. 

Similarly to informal mathematics, formal mathematical proofs necessarily draw on a large body of prior formalized results. The constraints of history loom even larger because these formal proofs must account for every inference from the axioms to the theorem in question: $\delta$ must compile in the context of all the prior results in the library---in the case of \texttt{Lean}, for instance, the \texttt{mathlib4} ecosystem.

Proofs and declarations are added to these libraries by members of a community defined by its commitment to that particular ecosystem. Unlike the ecosystem for $\pi$, which, due to its informal nature, can be pluralistic and tolerant in its conventions, formal libraries for $\delta$ are rigid dependency structures in which every definition must be precisely compatible with everything else.

The infrastructure required to prove strong results emerges as a by-product of the particular formalization projects the community undertakes. While some are the product of individual preference and whimsy, many are high-profile endeavours to demonstrate the promises of formalization to other mathematicians (as in, for example, the Liquid Tensor Experiment), funding agencies, or corporations who might hope to use the technology for commercial ends. Regardless of the motive, however, the result is that these early formalizations can become generatively entrenched and costly to change \cite{wimsatt2007re}. This kind of path dependence is already visible in formal proofs, whose dependency structures mirror those found in large-scale software projects such as the Linux kernel~\cite{viteri2022epistemic}. 

As the dependency graph grows around these earlier choices, later formalizations necessarily inherit them, regardless of how well they represent the underlying proof-ideas within the current proof being formalized. As a result, a mathematician working from an informal proof $\pi$ may be forced to use library definitions that fail to match the conceptual apparatus of $\pi$ even though they succeed in securing $\theta$, and the resulting formal derivation $\delta$ may only track $\pi$ weakly, if at all. Correspondence is hybridized, and the $\pi$-$\delta$ relationship is now mediated and altered by the conventions and constraints of the library---choices the formalizer inherits but did not make. Mathematics has always had a history, and a path dependence; what is new is its encounter with a novel, independent, and potentially less evolvable tradition. 

\section{Routes to Correspondence}
\label{routes}

How is correspondence established? We see at least three independent modes by which mathematicians working to formalize informal proofs can verify, in practice, that an adequate representation and tracking has been achieved. Following Wimsatt \cite{wimsatt2007re}, we take the independence of these modes to be epistemically significant. Convergence across modes that are individually vulnerable to different sources of failure provides a robustness that no single mode can deliver on its own. 

The first mode is \textit{structural}: mathematicians can verify that a $\pi$-$\delta$ pair share certain relevant ``structures''. The second is \textit{causal}: mathematicians can verify a causal pathway from the informal proof to the formal counterpart. The third mode is \textit{explanatory}: mathematicians can verify that the informal proof helps explain why the formal proof type-checks, or that the formal proof can help explain why the informal one works.

In principle, one establishes correspondence-by-structure by checking whether the objects and inferential moves in $\delta$ mirror those in $\pi$. In the simple {\tt primes\_unbounded} example above, one might notice the parallel between the informal statement to consider an $m$ equal to $n!+1$ and the opening definition of an ``{\tt m}'' as ``{\tt Nat.factorial n + 1}''---along with some additional stress-testing of {\tt Nat.factorial} and {\tt +} to ensure that these terms are doing what we expect. As we note above, this structural mode is likely to become increasingly fraught as the gap between the new mathematics we want to verify, and the mathematics the libraries were first built to verify, widens. If correspondence is to be established at all, it must draw on resources beyond structural comparison alone.

When we see formalization as translation, the \textit{causal} mode becomes relevant. In this mode, one sees how a \textit{particular step} in $\pi$, say, leads one to program a \textit{particular step} (and not some other, equally valid, one) in $\delta$, in the same way one ``follows'' an argument. Causal relationships need not preserve structure. If $\pi$ is written from a set theoretic or classical-logic perspective, but $\delta$ uses type theory and constructive logic, some steps in the mathematician's translation may bear only a weak structural relationship. It is also, at least at this stage, unclear what aspects of the causal relationship matter---particularly when A.I. tools are in play. The causal mode delivers \textit{provenance} and distinguishes a derivation that was produced by following this particular proof from one that merely happens to prove the same theorem.

\textit{Explanatory} relationships, finally, are not only a motivation for undertaking formalization in the first place \cite{avigad2024mathematics}, but can also be particularly powerful if, say, $\delta$ corresponds closely enough to $\pi$ that we can see how it fills in the gaps, or when, conversely, $\pi$ helps us make sense of a particularly opaque section of $\delta$. If a formalization deepens our understanding of an original informal proof---perhaps by drawing attention to less-noticed aspects of its underlying argument---we are likely correct in increasing our confidence that both representation and tracking have been achieved. This is the opposite of the old saw ``garbage in, garbage out'': when a formalization provides us with unexpected epistemic benefits for the original informal proof, this clues us in to the possibility that both representation and tracking were successful.

Both the causal and explanatory modes arise in Avigad's response to the controversy surrounding \emph{Math, Inc.}'s AI automated completion of the sphere packing proof. 

\begin{quote}
    ``[T]he formalization, on its own, is close to worthless, since the correctness of Viazovska’s result was never in doubt. The participants embraced the project, rather, as a way of revisiting those results and better understanding them, and of building libraries and infrastructure to support future work'' \cite[p.3]{avigad_sphere}.
\end{quote}

For the formalization to be worthwhile, in other words, it has to have been involved in a causal process (the revisiting of Viazovska's results through the activity of new human interpreters) and, one hopes, participate in an epistemic goal of increasing our mathematical understanding \cite{thurston2006proof}.

\section{Conclusion}
\label{ai}

The claim of this paper is that the promise of the Standard View---to explain the deductive certainty of our informal mathematical achievements---is a promise indefinitely postponed. Behind its deductive front is a much richer notion of correspondence that becomes tractable for finite minds by the convergence of multiple defeasible modes. Following Wimsatt, we might say that formalization advances mathematics precisely because these modes can fail independently. Formalization is not a deductive bridge from informal proof to formal certainty, but a defeasible, empirical achievement of correspondence supported by multiple independent evidential routes. 

The rapid advance of AI-guided mathematics brings these questions more urgently than ever to the fore \cite{duede2024apriori}. Once causal provenance and explanatory transparency are weakened, the convergence needed to sustain claims of thick correspondence is weakened as well.

\section*{Acknowledgments}

The Authors wish to thank Jeremy Avigad, David Barack, Daniel C. Friedman, Sidharth Hariharan, and Chase Norman for their helpful feedback on an earlier draft of this manuscript. This work was supported by Grant 63750, ``Explaining Universal Truths'', to SD from the John Templeton Foundation, and Grant G-2025-79234, ``AI \& Evolving Disciplinary Norms of Inquiry'', to ED from the Alfred P. Sloan Foundation.

\bibliographystyle{alpha}
\bibliography{bibliography}

\end{document}